\documentclass[12pt]{amsart}

\usepackage{amssymb}
\usepackage{tikz}
\usepackage{tikz-cd}
\usepackage{fullpage}
\usepackage{microtype}
\usepackage{hyperref}

\usepackage{multicol} 
\usepackage{lineno}   

\newcommand{\lean}[1]{\lstinline[language=lean]{#1}}

\definecolor{keywordcolor}{rgb}{0.7, 0.1, 0.1}   
\definecolor{tacticcolor}{rgb}{0.0, 0.1, 0.3}    
\definecolor{commentcolor}{rgb}{0.4, 0.4, 0.4}   
\definecolor{stringcolor}{rgb}{0.5, 0.3, 0.2}    
\definecolor{symbolcolor}{rgb}{0.1, 0.2, 0.7}    
\definecolor{sortcolor}{rgb}{0.1, 0.5, 0.1}      
\definecolor{attributecolor}{rgb}{0.7, 0.1, 0.1} 
\definecolor{errorcolor}{rgb}{1, 0, 0}           

\usepackage{listings}

\usepackage{scalefnt}

\newcommand{\sorry}[0]{\lean{sorry}}

\newcommand{\Cond}[0]{\operatorname{Cond}}
\newcommand{\Ab}[0]{\mathrm{Ab}}
\newcommand{\Mcal}[0]{\mathcal{M}}
\newcommand{\Hom}[0]{\operatorname{Hom}}
\newcommand{\Ext}[0]{\operatorname{Ext}}
\newcommand{\Prof}[0]{\mathrm{Profinite}}
\newcommand{\Zbb}[0]{\mathbb{Z}}
\newcommand{\Qbb}[0]{\mathbb{Q}}
\newcommand{\Rbb}[0]{\mathbb{R}}

\newtheorem*{theorem*}{Theorem}

\newtheorem{maintheorem}{Theorem}

\title{Abstraction boundaries and spec driven development in pure mathematics}
\author{Johan Commelin}
\address{Albert--Ludwigs-Universit\"at Freiburg, Mathematisches Institut, Ernst-Zermelo-Stra\ss{}e~1, 79104 Freiburg im Breisgau, Deutschland}
\email{jmc@math.uni-freiburg.de}
\thanks{JC is supported by the Deutsche Forschungsgemeinschaft (DFG) grant CO~2587/1-1.}
\author{Adam Topaz}
\address{Mathematical and Statistical Sciences, University of Alberta, 632 Central Academic Building, Edmonton AB T6G 2G1, Canada}
\email{topaz@ualberta.ca}
\thanks{AT is supported by NSERC discovery grant RGPIN-2019-04762.}
\date{\today}

\begin{document}
\maketitle

\begin{abstract}
  In this article we discuss how \emph{abstraction boundaries} can help tame complexity in mathematical research, with the help of an interactive theorem prover.
  While many of the ideas we present here have been used implicitly by mathematicians for some time, we argue that the use of an interactive theorem prover introduces additional qualitative benefits in the implementation of these ideas.
\end{abstract}

\section{Introduction}

Modern research in pure mathematics has a clear tendency toward increasing complexity.
New striking mathematical results may involve complex proof techniques, deep mastery of a subfield within mathematics, or nontrivial input from several areas of mathematics.
All of these play a role in increasing the inherent complexity of a piece of modern mathematics.

In many respects, such increasing complexity is an indication of \emph{progress} in pure mathematics.
However, with such complexity comes a significant increase in \emph{cognitive load} for both readers and authors.
It is now routine to see significant new papers in pure mathematics with one hundred pages or more.
Similarly, the refereeing process of a complex mathematical result regularly takes multiple years.
This state of affairs is clearly not sustainable.

Mathematicians have long been using the concept of abstraction in order to tame such complexities.
For example, they often break down proofs into more manageable steps, such as lemmas, propositions, etc.
They also frequently introduce definitions that capture a curated selection of properties satisfied by the objects they wish to study.
More generally, it is often helpful to create ``abstractions boundaries'', which separate the use of a mathematical object from its actual implementation.
Staying within the boundaries of a given abstraction can help reduce cognitive load to some extent, but inherent complexities may nevertheless remain.

In this article we discuss how the use of abstraction boundaries, in conjunction with an \emph{interactive theorem prover} (ITP), can facilitate complex interdisciplinary mathematical collaboration.
We argue that the implementation of these ideas within an ITP can have significant qualitative benefits in further reducing cognitive load for both authors and readers.
The application of this methodology was a key ingredient in the success of the \emph{liquid tensor experiment} (LTE), and we use examples from this project to illustrate our primary points.

Here is a brief outline of the paper.
In \S2, we review some aspects of condensed mathematics and LTE, focusing primarily on the topics which are relevant for this paper.
In \S3, we give a more precise description of what we mean by an \emph{abstraction boundary} along with a discussion on sources of complexity within mathematics.
In \S4, we highlight some of the qualitative benefits offered by the use of abstraction boundaries alongside an ITP, while illustrating \emph{spec driven development} with several examples related to LTE.
The issue of aligning formal and informal mathematics using abductive reasoning is addressed in \S5.
Finally, we summarize and conclude our discussion in \S6.


\subsection*{Acknowledgements}
We would like to warmly thank all of the contributors to the liquid tensor experiment for the amazing collaboration and excellent discussions.
We thank the Lean community at large for providing the technological and social environments which made this work possible.
Many thanks to the referees of this article and the editors of this volume for their comments which helped improve this paper in various ways.
Finally, we would like to thank A.~Venkatesh for proposing such an exciting topic for his Fields Medal Symposium.

\section{An overview of the liquid tensor experiment}

The \emph{liquid tensor experiment} (LTE) began with a challenge posed by P.~Scholze~\cite{zbMATH07566886}, whose ultimate goal was to formally verify the proof of the \emph{main theorem of liquid vector spaces}:
\begin{maintheorem}[Clausen--Scholze]\label{theorem:challenge}
  Let $0 < p' < p \le 1$ be real numbers and let $S$ be a profinite set.
  Let $V$ be a $p$-Banach space, considered as a condensed abelian group.
  Write $\Mcal_{p'}(S)$ for the space of real-valued $p'$-measures on $S$, also considered as a condensed abelian group.
  Then $\Ext^{i}_{\Cond(\Ab)}(\Mcal_{p'}(S),V) = 0$ for all $i \geq 1$.
\end{maintheorem}

This theorem is of utmost foundational importance; it is the technical core of the assertion that the real numbers, endowed with the collection of so-called \emph{$p$-liquid vector spaces}, forms an \emph{analytic ring} in the sense of Clausen--Scholze~\cite{condensedpdf,analyticpdf,complexpdf}.
In~\cite{zbMATH07566886}, Scholze mentions that this might be his ``most important theorem to date.''
In this section, we will only give a brief overview and motivation for the objects involved in this theorem along with a general summary of the experiment itself, as will be relevant for the remainder of this paper.
We refer the reader to~\cite{zbMATH07566886} and to~\cite[Theorem 9.1]{analyticpdf} for more details.

\subsection{Condensed and liquid mathematics}

Condensed mathematics is a new foundational system which replaces topological spaces with objects called \emph{condensed sets}, which are set-valued sheaves on the category $\Prof$ of profinite sets with respect to the coherent Grothendieck topology.
The power of this new theory shines particularly in algebraic contexts, whereas the classical approaches using point-set topology have several defects.
For example, the category of topological abelian groups fails (quite badly) to be an abelian category, while the category of condensed abelian groups, defined either as sheaves of abelian groups on $\Prof$ or equivalently as the category of abelian group objects in condensed sets, is an exceptionally nice abelian category.
Condensed mathematics and its applications are still under development by Clausen--Scholze~\cite{condensedpdf,analyticpdf,complexpdf} et al., while a similar approach, albeit with a different overall focus, was introduced by Barwick--Haine in~\cite{barwick2019pyknotic} under the name \emph{pyknotic sets}.%
\footnote{The actual definition of a condensed set involves cardinality bounds which we have omitted from the discussion, while the only distinction between condensed and pyknotic objects is in such set-theoretic issues.}

Any topological space (resp.~group, ring, module, etc.) can be considered as a condensed set (resp.~group, ring, module, etc.) via a Yoneda-like construction.
One may thus speak about condensed modules over condensed rings, which are the condensed incarnation of topological modules over topological rings.
But in order to capture the building blocks of \emph{analytic geometry}, one needs a notion of \emph{completeness} for such modules that behaves well algebraically.
This is exactly what the notion of an \emph{analytic ring} entails: it is a condensed ring with a notion of ``completeness'' for condensed modules which behaves well in a suitable \emph{derived} sense.

Examples of analytic rings include $\Zbb^{\blacksquare}$, whose underlying condensed ring is $\Zbb$ endowed the discrete topology, where the notion of completeness for modules is modeled on $\Zbb[S]^{\blacksquare} := \varprojlim_{i} \Zbb[S_{i}]$ for profinite sets $S = \varprojlim_{i} S_{i}$ expressed as cofiltered limits of finite sets $S_{i}$.
A similar construction makes $\Zbb_{p}$, with the $p$-adic topology, into an analytic ring denoted $\Zbb_{p}^{\blacksquare}$.
The field $\Qbb_{p}$ can be considered as an analytic ring as well, essentially by base-changing from $\Zbb_{p}$.

The prototypical complete $\Qbb_{p}$-module in the last case turns out to be $\Mcal(S,\Qbb_{p})$, the space of (bounded) $\Qbb_{p}$-valued measures on a profinite set $S$.
However, the analogous assertion \emph{fails} if one replaces $\Qbb_{p}$ with $\Rbb$.
The theory of \emph{$p$-liquid vector spaces} solves this issue, by replacing the prototype $\Mcal(S,\Rbb)$ by a \emph{smaller} space $\Mcal_{< p}(S,\Rbb)$.

Theorem~\ref{theorem:challenge} eventually translates to the assertion that $\Rbb$, endowed with the collection of such $p$-liquid vector spaces, forms an analytic ring, thus providing a family of analytic rings, parameterized by $0 < p \le 1$, whose underlying condensed ring is $\Rbb$.
This result shows real and complex geometry can be studied in the context of analytic spaces alongside objects from algebraic and nonarchimedean geometry, see~\cite{analyticpdf,complexpdf}.
At the same time, classical objects from functional analysis, such as Banach spaces (or more generally, $p$-Banach spaces for $0 < p \le 1$), can fit in this framework as well, and thus many facets of functional analysis can now be studied \emph{algebraically} with the necessary analysis essentially encapsulated in Theorem~\ref{theorem:challenge}.
This work has already resulted in generalizations and new \emph{algebraic} proofs of classical results from algebraic geometry which were originally analytic in nature, such as Serre's GAGA.
Many additional applications will surely come in the future.

At this point, it should be quite clear that the proof of Theorem~\ref{theorem:challenge} involves significant amounts of category theory, topos theory, and homological algebra, as well as analysis and topology.
Surprisingly, a key step in the proof (see~\cite{386796,Scholze-half-a-year}) brings combinatorics and (discrete) convex geometry into the picture as well.

\subsection{Verifying the theorem on a computer}
The liquid tensor experiment (LTE) was completed on July 14, 2022~\cite{LTE-completion} with a complete formal verification of Theorem~\ref{theorem:challenge} using the Lean proof assistant~\cite{10.1007/978-3-319-21401-6_26} and its mathematics library \texttt{Mathlib}~\cite{10.1145/3372885.3373824}.
The formalization project, which was led by the authors of the present paper, was a collaboration involving about a dozen mathematicians and several computer scientists.
A detailed exposition of the project itself will appear elsewhere.

To conclude this section, we would like to highlight again that the proof of Theorem~\ref{theorem:challenge} involves many disparate subfields of pure mathematics, and at the same time, many relatively complicated mathematical objects.
This significant source of complexity in the proof would have posed serious challenges in a traditional refereeing process.
Indeed, this is one of the main reasons Scholze sought a formalization of the proof of Theorem~\ref{theorem:challenge} as explained in~\cite{zbMATH07566886}.
In~\cite{Scholze-half-a-year}, he writes:
\begin{quotation}
  The Lean Proof Assistant was really that: An assistant in navigating through the thick jungle that this proof is.
  Really, one key problem I had when I was trying to find this proof was that I was essentially unable to keep all the objects in my “RAM”, and I think the same problem occurs when trying to read the proof.
\end{quotation}
The Lean ITP helped wrangle this complex proof, and, as we explain below, abstraction boundaries played a key role in the process.


\section{Managing complexity with abstraction boundaries}\label{section:sources-of-complexity}


Complexity in mathematical research projects can come from various sources.
We distinguish between two kinds.
The first is \emph{inherent complexity}, by which we mean the raw complexity of the Platonic ideal of a piece of mathematics.
While this notion is hard to quantify, an approximate measure of inherent complexity could be the amount of time it takes a generic mathematician to understand the proof or definition in question starting from scratch.

The second form of complexity we call \emph{accidental complexity}.
This encompasses all sources of complexity that are not essential to the mathematics
but rather imperfections that could be removed given enough time and energy.

Both inherent and accidental complexity contribute to the total complexity of a piece of mathematics.
It is this total complexity which can overwhelm both authors and readers.

\subsection{Inherent complexity}
The main source of inherent complexity for a piece of mathematics comes from the scope of its prerequisites.
This can be measured along the two axes of \emph{breadth} and \emph{depth}:
a piece of mathematics may combine ingredients from many different areas of mathematics,
or it may be a pinnacle within some subject, or a combination of both.

Note that the complexity of a theorem statement can be unrelated to the complexity of its proof.
Fermat's Last Theorem is an excellent example in this regard.
The statement can be understood by anyone who understands the usual operations on natural numbers,
but understanding of the proof requires mastery of several subfields of mathematics.
Conversely, an elementary property of a complicated definition may have a simple proof.
In this case, the \emph{statement} has a high inherent complexity, but the proof does not.

Moreover, two different approaches to a proof of a given theorem may differ in their inherent complexities, even if both are written in an ideal manner.
Thus, inherent complexity does not measure the minimal complexity of any possible proof or argument,
but rather the complexity of a particular argument or proof strategy.
An excellent proof strategy that is executed poorly may have low inherent complexity and high accidental complexity (see below).
Vice versa, a poorly chosen proof strategy that is executed perfectly may have high inherent complexity and low accidental complexity.

Complexity may also arise from the applications of proof \emph{techniques} as opposed to applications of theorems.
By this we are referring to proofs that apply variations of a common technique, but where it is unreasonable to capture all applicable variations in single general theorem statement.
This happens frequently in parts of analysis, probability theory, and combinatorics, to name some examples.
A specific example, provided to us by Peter Nelson, is ``connectivity reduction,'' where one proves theorems about combinatorial objects such as graphs, matroids, hypergraphs, etc., by reducing to certain ``connected'' objects.
There are many variations in the objects involved, the notion of connectivity, the way connected objected are glued together, etc.
While such connectivity reductions are common, it does not seem reasonable to unify them in a precise theorem statement.

While it may appear as if this is more a form of accidental complexity, we argue that these proof techniques should be rather viewed as inherent complexity of the subfield of mathematics in which the proof takes place.
Although it might be possible to rephrase such a proof technique as an actual theorem for the purpose of a given proof or mathematical argument,
the theorem will not be sufficiently broadly applicable: only the proof strategy generalizes.

\subsection{Accidental complexity}
We illustrate accidental complexity by some well-known examples that occur as part of the reading and writing of mathematical texts.

\emph{Referring to the proof of a lemma, instead of its statement.}
This may occur in the rather common situation where
the proof of Lemma~X.Y.Z in reference~[AB] actually proves something stronger than what the statement of Lemma~X.Y.Z claims.
In this case we might write something like
``We now see that our quasi-perfect gadget~$G$ is strongly separated, by the proof of [Lemma~X.Y.Z, AB].
Indeed, it is not hard to check that the proof shows $G$ is not only separated, but even \emph{strongly} separated''.

\emph{Using modifications of source material.}
Sometimes in the mathematical writing process,
we need a verbatim copy of some source material modulo a handful of completely algorithmic modifications.
A mild example of this might be:
``In the following proof, we will rely on a slightly modified version of the material in section~5 of~[CD],
but with the condition that the field $k$ has characteristic~$0$ replaced by the condition that $k$ is perfect.''
This can often occur in combination with the preceding point, where a modification of a proof of some result applies while the statement itself does not.

\emph{Unwritten assumptions.}
Sometimes statements of theorems or lemmas can omit assumptions.
In the good cases, these assumptions were listed at the beginning of the section
(``in this section we assume that $G$ is infinite'')
or at the beginning of the paper in a section on notation.
But there might also be assumptions that are common practice within certain communities
(``all rings are commutative'')
which might be hard to discover for the uninitiated reader.

\emph{Unclear definitions or terminology.}
A related point is that it may be hard to understand definitions or terminology.
For example, definitions may change over time.
In the early literature on scheme theory, all schemes were separated.
In Whitney's excellent book on Geometric Integration Theory~\cite[p.15]{whitney},
a function is called smooth if it is $C^1$.
Such changes occur naturally over time,
but can lead to significant confusion when studying older works,
and even to a second degree when contemporary sources quote older works without properly warning about the change in terminology.

\emph{Keeping track of side conditions.}
When a lemma is applied, we verify the main assumptions of the lemma.
But it would be burdensome and distracting to verify all the ``trivial'' side conditions
(that some number is positive, or less than~$\epsilon$, or that some map is continuous, etc).
At the same time, requiring the reader to keep track of these conditions
implicitly increases the cognitive load.

\emph{Editing a paper.}
Making changes to a paper in its final stages can be tricky,
because of all the possible ramifications downstream.
Moving a lemma from section~2 to section~3 might require changing the statement,
because the standing assumptions on~$X$ differ between the two sections.
And yet it is all too easy to make such changes in a careless manner,
possibly even shortly before publication and after a refereeing process took place.

All of the examples above, among others, are well-known pitfalls in mathematical writing, and are often warned against, see \cite{halmos1970write, knuth1989mathematical, serre2007badly}.
Nevertheless these issues are still commonplace in the mathematical literature, increasing the accidental complexity of the work.

\subsection{What is an abstraction boundary?}
An \emph{abstraction boundary} is a formal separation between the implementation details of a concept from its extrinsic properties and behaviour.
This involves introducing a \emph{specification} (or \emph{interface}) which describes how the object interacts with the outside environment.
By using such specifications, one can work with the concepts they capture without relying on any actual implementation details.

This concept is prevalent in software engineering, with prominent examples including \texttt{C} header files, public methods in object-oriented programming, and typeclasses in functional programming.
In mathematics, abstraction boundaries play a fundamental role as well,
with the mathematical concept of a \emph{definition} being the principal example.
For example, the definition of a ``group'' allows us to develop group theory, and apply results from it in the most varying circumstances,
instead of only working with particular examples such as permutation groups and matrix groups.
In this context, the specification is provided by the group axioms.
Another example of abstraction boundaries in mathematics includes ``black box'' theorems, such as Cauchy's residue theorem, the existence of a N\'eron model, the law of large numbers, Theorem~\ref{theorem:challenge}, etc.
Such theorems can be effectively applied in particular situations without requiring an understanding of the argument.

\subsection{Abstraction boundaries in interactive theorem provers}
Interactive theorem provers (ITPs) enforce certain abstraction boundaries rigidly.
While this can be restraining at times, it also relieves the user from cognitive load arising from several sources of complexity.

ITPs differ in the features that they offer,
and certainly in the implementation details of those features.
But most of this section applies to all modern ITPs,
even though we use the Lean theorem prover as a running example.
Lean is the system that we used for LTE,
and it is the ITP that we are most familiar with.

The most well-known feature offered by ITPs is \emph{proof checking}, often as an instance of \emph{type checking}.
This means that the ITP ensures that lemmas and theorems are always applied correctly,
and that all side conditions are satisfied.
Certain side conditions that are amenable to automation, such as continuity proofs, functoriality proofs, etc., can often be handled in a way that is transparent to the user.

Another feature that most ITPs offer is \emph{proof irrelevance}.%
\footnote{ITPs built on a foundation of homotopy-type-theoretic nature will often not have proof irrelevance for \emph{explicit reasons} arising from the foundational theory.
The particular reasons for this choice are not relevant to this paper,
and a discussion of homotopy type theory would lead us too far astray.}
This notion can be seen as an analogue of the common mathematical practice of treating certain theorem statements as black boxes.
Essentially, this means that after an ITP finishes checking the proof of a lemma or theorem,
it remembers only one bit of information: that the proof is valid.
Any constructions performed within the proof become inaccessible as soon as the proof is complete.
If it is necessary to refer to such a construction later on,
it must be factored out of the proof into a standalone definition.
In doing so, the construction becomes an object in its own right,
that can be named, referred to, and reasoned about.
Conversely, if a construction is local to a proof,
it implicitly comes with the strict promise that it will not be used outside of the proof.

\section{Spec driven development}\label{section:spec-driven-dev}

\emph{Spec driven development} is an approach for ``doing mathematics'' which aims to both use and capture the abstraction boundaries of the mathematical objects involved.
For our purposes, by ``doing mathematics'' we mean the process of making mathematical objects such as definitions, constructions, etc., and proving theorems about these objects.
As we explain below, the use of an ITP in conjunction with spec driven development provides substantial benefits in taming the various sources of complexity, particularly in complex and collaborative projects.

While this approach arose naturally during the LTE, it is quite commonplace, although \emph{implicit}, in usual mathematical practice.
For example, an informal proof may assume some properties up front with the promise of a justification in a later section.
The main difference here is that organizing a proof in this way usually increases the accidental complexity and hence the cognitive load.
For example, the reader may feel inclined to check that there is no circular reasoning in the final argument.
The use of an ITP completely eliminates this source of complexity, and, as explained below, even \emph{encourages} such uses of mathematical debt.

Spec driven development can be summarized as the following recursive process:
\begin{enumerate}
  \item Isolate a desired mathematical definition, construction, theorem, etc., as a ``target''.
  \item When appropriate, isolate an initial \emph{specification} (``spec'' for short) of the target.
  \item Break down the object and its specification into parts with lower complexity.
  \item Repeat the above, with each new part acting as the target.
\end{enumerate}
Furthermore, at each step, definitions and/or theorems, including the targets and/or their specs, may be \emph{refactored} as necessary.
Refactoring refers to the process of restructuring the implementation of an existing object without changing its external behaviour.
As long as these implementation details and external properties of the object are properly separated by an abstraction boundary, other parts of the environment will not be affected by such a restructuring.

Spec driven development is also deeply intertwined with the notion of an abstraction boundary discussed in the previous section.
In broad terms, the relationship between the two can be broken up into two interconnected categories:
\begin{enumerate}
  \item Spec driven development uses abstraction boundaries in order to control complexity.
  \item New abstraction boundaries arise naturally in the process of spec driven development.
\end{enumerate}
Of course, these two categories are far from disjoint.
For instance, as item (2) produces abstraction boundaries, they naturally feed into item (1).
At the same time there are situations where spec driven development can straddle both categories.

In the following subsections we discuss the idea behind spec driven development in more detail by going through a series of examples, based on constructions that were necessary in the LTE.
We will present a few snippets of Lean code, although experience with the language or any other ITP will not be necessary.
The precise meaning of the code we present and the mathematics involved will not be important, and are only meant to illustrate the implementation of spec driven development in an ITP.
The only important concept needed from Lean's syntax is that of a \sorry; this is a keyword which should be thought of as a \emph{placeholder} that needs to be filled in later.
A \sorry~can be used as a placeholder for both \emph{data} and \emph{proofs}, which is crucial for this approach.\footnote{In type-theoretic terms, a \sorry~can be used as a placeholder for an arbitrary \emph{term} of an arbitrary type, including propositions.}

\subsection{Using and extending abstraction boundaries}\label{subsection:abelian-cat-example}
The simplest use-case of spec driven development is in using previously established abstraction boundaries in order to reduce complexity for the overall objective at hand.
As an example, consider the assertion that the category of condensed abelian groups is an abelian category, a property that is used in numerous places in LTE.
In Lean, this might be initially formulated as follows:
\begin{lstlisting}
instance : abelian (Condensed Ab) := sorry
\end{lstlisting}
In this example, the \sorry~appearing should be considered as our initial \emph{target}: it is a unit of mathematical debt that must be accounted for at some point in the future.
At this point, all collaborators can immediately start using the assertion that the category \lean{Condensed Ab} is abelian, even before providing a proof of this fact.
For example, we can now speak about kernels and cokernels of morphisms, exact sequences, etc.
The abstraction boundary here is the notion of an \emph{abelian category}.

There are instances where we wish to extend the boundary in some way, and at the basic level, this can be accomplished by adding additional \emph{spec lemmas/theorems}.
For instance, in the example above we may want to use the fact that morphisms of condensed abelian groups are morphisms of presheaves on \lean{Profinite} taking values in \lean{Ab}, and that sums of such morphisms are computed ``object-wise''.
In practice, this assertion would be true by definition, but since the actual definition has yet to be provided, we cannot formally rely on this.
Instead, we add this property with a \emph{spec lemma}, which becomes an additional target that can again be used immediately:
\begin{lstlisting}
lemma val_app_add {F G : Condensed Ab} {X : Profiniteᵒᵖ}
  (f g : F ⟶ G) : (f + g).val.app X = f.val.app X + g.val.app X :=
sorry
\end{lstlisting}
In the code above, if
\lean{f}
is a morphism of condensed abelian groups, then \lean{f.val.app X} denotes the corresponding morphism on sections over the profinite set \lean{X}.
Note that if we \emph{use} the targets above at this point, then we are \emph{forced} to remain within their associated abstraction boundary.
We can't rely on an implementation detail because there is no implementation yet!

In order to start repaying the existing debt, we can break down each target and its associated spec into smaller parts, and repeat the process until the individual targets can be handled directly.
The initial target and its spec may thus become a collection of several parts, each of which becomes an additional target.
In our example, recall that an abelian category is a preadditive category satisfying some additional conditions.
After we do this for the abelian category instance in the example above, we are left with the following which summarizes all the current targets:
\begin{lstlisting}
instance : preadditive (Condensed Ab) := sorry

instance : abelian (Condensed Ab) :=
{ normal_mono_of_mono := sorry,
  normal_epi_of_epi := sorry,
  has_finite_products := sorry,
  has_kernels := sorry,
  has_cokernels := sorry,
  ..(infer_instance : preadditive (Condensed Ab)) }

lemma add_val_app {F G : Condensed Ab} {X : Profiniteᵒᵖ}
  (f g : F ⟶ G) : (f + g).val.app X = f.val.app X + g.val.app X :=
sorry
\end{lstlisting}

At each step, a contributor may fill in some target(s) by replacing the \sorry~with an actual construction/proof.
For example, one contributor may construct the preadditive structure on \lean{Condensed Ab} in such a way that the spec lemma \lean{add_val_app} above is true by definition, resulting in a reduction of debt:

\begin{lstlisting}
instance : preadditive (Condensed Ab) := /- actual proof omitted -/

lemma add_val_app {F G : Condensed Ab} {X : Profiniteᵒᵖ}
  (f g : F ⟶ G) : (f + g).val.app X = f.val.app X + g.val.app X :=
rfl
\end{lstlisting}

At this point, the remaining debt in our running example consists of \lean{normal_mono_of_mono}, \lean{normal_epi_of_epi}, etc., which appear as components of the proof that \lean{Condensed Ab} is an abelian category.
These targets may again be handled directly, or further broken down as necessary.

\subsection{Related targets and specs}\label{subsection:ext-example}

In practice, the targets appearing in spec driven development may depend in nontrivial ways on others.
To illustrate this, let's consider a more interesting example that came up during the LTE.
In order to state Theorem~\ref{theorem:challenge} formally, we must be able to talk about $\Ext$-groups, and a natural context for this construction is using an abelian category with enough projectives.
For each natural number $n$, $\Ext^{n}(-,-)$ is a bifunctor taking values in abelian groups, which is contravariant in the first variable.
We may thus consider the following as our target:

\begin{lstlisting}
-- The arrow `⥤` is used to denote functors.
def Ext {A : Type*} [category A] [abelian A] [enough_projectives A] (n : ℕ) :
  Aᵒᵖ ⥤ A ⥤ Ab :=
sorry
\end{lstlisting}

In this case, it makes sense to immediately introduce specifications describing the intended behavior of this definition.
For example, we may want to ensure that $\Ext^{0}$ is naturally isomorphic to $\Hom$, and to speak about long exact sequences of $\Ext$-groups, say in the first variable.
This can all be accomplished by adding the following:
\begin{lstlisting}
-- `preadditive_yoneda.flip` is the bifunctor `X ↦ (Y  ↦ Hom(X,Y))`,
-- where `Hom(X,Y)` is considered as an abelian group.
def Ext_zero_iso_Hom {A : Type*} [category A]
  [abelian A] [enough_projectives A] :
    (Ext 0 : Aᵒᵖ ⥤ A ⥤ Ab) ≅ preadditive_yoneda.flip :=
sorry

def d {A : Type*} [category A] [abelian A]
  [enough_projectives A] (n : ℕ) (X Y : A) :
    ((Ext n).obj (op X)).obj Y ⟶ ((Ext (n+1)).obj (op X)).obj Y :=
sorry

lemma Ext_LES {A : Type*} [category A]
    [abelian A] [enough_projectives A] (n : ℕ)
    (X₁ X₂ X₃ Y : A) (f : X₁ ⟶ X₂) (g : X₂ ⟶ X₃)
    (ses : short_exact f g) :
  exact_seq Ab
    [((Ext n).map g.op).app Y, ((Ext n).map f.op).app Y, d n X₁ Y] :=
sorry
\end{lstlisting}
The definitions \lean{Ext_zero_iso_Hom}, \lean{d} and spec lemma \lean{Ext_LES} are additional targets, while \lean{Ext_LES} relates to \emph{all three} targets \lean{Ext}, \lean{Ext_zero_iso_Hom} and \lean{d}.

As additional targets and/or preexisting definitions/theorems are introduced into the context, the situation can quickly become unwieldy.
It is at this stage where the use of an ITP really shines.
Essentially, the ITP keeps track of all the connections between the targets and existing definitions/theorems, while ensuring consistency \emph{at all times.}
One enormous qualitative benefit of using an ITP with this approach is that a contributor working on a given target only needs to keep track of the immediately relevant portions of this dependency graph of targets/definitions/theorems, thereby substantially reducing cognitive load.
This approach also facilitates collaboration, as it is no longer necessary to remember (or even understand) the context of the project as a whole when working on individual targets.
In collaborative projects involving various mathematical subfields, each contributor may thus choose to only work on targets within their own area(s) of expertise.

\subsection{Refactoring}\label{subsection:refactoring}
Another major benefit of using abstraction boundaries, which is highlighted in the process of spec driven development, is that \emph{refactoring} becomes easier in most cases.
For instance, in the example from \S\ref{subsection:abelian-cat-example}, the actual \emph{implementation} of the fact that \lean{Condensed Ab} is an abelian category is irrelevant, as long as one only needs results which hold true in an arbitrary abelian category.
Changing the implementation details will therefore not interfere with the validity of such results.

Continuing with the example of $\Ext$ groups from \S\ref{subsection:ext-example}, in LTE there was a point where we had a working definition of \lean{Ext}, similar to the one described in \S\ref{subsection:ext-example}, which had to be refactored to accommodate $\Ext$-groups of (bounded above) \emph{complexes} as opposed to just objects.
In other words, at one point we introduced targets resembling the following:
\begin{lstlisting}
def bounded_above_complex.Ext {A : Type*} [category A]
    [abelian A] [enough_projectives A] (n : ℤ) :
  (bounded_above_complex A)ᵒᵖ ⥤ (bounded_above_complex A) ⥤ Ab :=
sorry

def bounded_above_complex.single_zero {A : Type*} [category A] [abelian A] :
  A ⥤ bounded_above_complex A :=
sorry
\end{lstlisting}
and \emph{redefined} the \lean{Ext} from \S\ref{subsection:ext-example} using these roughly as follows:
\begin{lstlisting}
def Ext {A : Type*} [category A]
    [abelian A] [enough_projectives A] (n : ℕ) :
  Aᵒᵖ ⥤ A ⥤ Ab :=
(bounded_above_complex.single_zero A).op ⋙
  (bounded_above_complex.single_zero A ⋙
    (bounded_homotopy_category.Ext n).flip).flip
\end{lstlisting}
By relying on the given abstraction boundaries, as one is essentially \emph{required} to do in the context of spec driven development, the process of resolving errors arising as a consequence of such a refactor becomes significantly easier.
While this change in \lean{Ext} may indeed cause certain other targets/specs/definitions/theorems to fail, this failure will be \emph{localized} to those constructions and proofs that depend on the  \emph{implementation}.
Crucially, \emph{the ITP will catch all such errors!}
Any results elsewhere in the project that rely on the abstraction boundary in question, such as proofs that only rely on the lemma \lean{Ext_LES} from \S\ref{subsection:ext-example}, should still work without any further modification.

\section{Aligning the formal with the informal}

Every mathematician has to make an ongoing effort to align
internal mental models of mathematics with the pen-and-paper representations found in the literature.
This alignment goes both ways:
mental models have to be updated and adjusted,
and expositions of mathematics can be streamlined and improved.
This alignment is an integral part of the creative process.
It happens at the blackboard, during seminar talks, at the writing table, or during a walk in the park.
Furthermore, it cannot be captured in formal rules:
one can not \emph{prove} that a mental model aligns well with a pen-and-paper representation of some piece of mathematics.

\subsection{The alignment problem}
With the use of ITPs, this alignment issue is extended in a new direction:
besides mental models and the pen-and-paper representations,
there are the additional \emph{digital} representations.
To illustrate the alignment problem, consider the statement of Theorem~\ref{theorem:challenge} as it was formalized in LTE:
\begin{lstlisting}
variables (p' p : ℝ≥0) [fact (0 < p')] [fact (p' < p)] [fact (p ≤ 1)]

theorem liquid_tensor_experiment (S : Profinite.{0}) (V : pBanach.{0} p) :
  ∀ i > 0, Ext i (ℳ_{p'} S) V ≅ 0 :=
/- proof omitted -/
\end{lstlisting}
At a superficial level, this code snippet shows a lot of similarity with the statement of Theorem~\ref{theorem:challenge}.
However, there is a possibility that the definition of the symbols \lean{Profinite}, \lean{pBanach}, \lean{Ext}, and so forth,
do not mean what they should mean according to established mathematical tradition.
In other words: do their formal definitions align with the pen-and-paper representations and/or with the mental model of anyone looking at the source code of LTE?
If the contributors to LTE were evil, the symbol \lean{Ext} might be defined as the zero functor, trivializing the statement in the code snippet above.
More innocently, there may be some subtle mistake in the definition of \lean{Ext} that implies the statement for completely uninteresting reasons.
Even this is quite rare in practice, as a formal proof written by an honest human mathematician still relies on the mathematician's mental models of the objects involved.
Nevertheless, the alignment issue still needs to be addressed in order to efficiently \emph{communicate} the results of a formalization to the greater mathematical community.

\subsection{Abductive reasoning}
Although there is no formal system for verifying alignment of different representations of mathematics,
we argue that \emph{abductive reasoning}, along with purposefully designed abstraction boundaries,
can be fruitfully applied to provide convincing evidence of alignment.
Abductive reasoning is the technical term for a mode of reasoning that
seeks to explain a set of facts or observations in the simplest way possible,
and is also called \emph{inference to the best explanation}.
Unlike with deductive reasoning,
an abduction is not a formal logical consequence of its premises.
The method is concisely captured by the colloquial ``Duck Test'':
\begin{quotation}
  If it looks like a duck, swims like a duck,\\
  and quacks like a duck, then it probably is a duck.
\end{quotation}
We provide a curated selection of examples and lemmas,
somewhat analogous to ``test suites'' in software developments,
that exhibit the ways in which a particular formal definition can be used.
This is the input that allows others to abductively conclude that these formal symbols align with their mental models of the corresponding mathematical notions.
In the terminology of the preceding section, this input would be a carefully crafted ``spec'' for the formal representation of a mathematical concept.

The crucial aspect of such a spec is that it can be much shorter and easier to read than the actual definitions themselves.
Indeed, typically the definitions will recursively unwind to a long list of prerequisite definitions.%
\footnote{The main result of LTE depends recursively on several thousands of definitions.
Part of these definitions are made in the project itself, but the majority is imported from Lean's mathematical library \texttt{mathlib}.}
Furthermore, definitions might be made using general constructions that the reader is unfamiliar with.

For LTE we provided such specs for all the objects occuring in the main theorem of liquid vector spaces.
We collected these specs in the subfolder \lean{examples/} of the project; see~\cite{LTE-examples} for a detailed discussion.
For example, we show that the real numbers in Lean are a conditionally complete linearly ordered field,
that $\text{Ext}^1(\Zbb/n\Zbb, \Zbb/n\Zbb) \cong \Zbb/n\Zbb$,
and that $p$-Banach spaces~$V$ can be given a norm that satisfies $\|\lambda v\| = |\lambda|^p\|v\|$ for all $\lambda \in \Rbb$ and $v \in V$.
Such examples are a \emph{reality check}, which collectively are meant to provide a large degree of confidence that the formal representations of mathematical concepts align with corresponding mental models and pen-and-paper representations.

Abductive reasoning is not unique to LTE or the use of ITPs, and many forms are often used in the mathematical practice at large.
For example, the following common methods of establishing and/or strengthening belief in a mathematical argument or proof can be considered as instances of abductive reasoning:
\begin{itemize}
	\item the ingredients and/or ideas seem strong enough to prove the claim;
	\item the claims hold up in well-chosen examples or numerical experiments;
	\item the proof is accepted by the community, including the experts in the field.
\end{itemize}
With such methods, a reader may obtain some level of trust in the validity of a piece of mathematics without necessarily digging into every detail of the exposition.
The reader must balance their desired level or trust with the amount of effort they are willing to exert.

The main differences between such forms of abductive reasoning and the strategy we put forth in LTE is in the scope to which it is applied.
While the three examples above focus on deciding whether or not to believe the \emph{proofs} of certain results, in the case of formalized mathematics, the focus shifts to merely the validity of the \emph{definitions} involved.
Thus, in formal projects such as LTE, the ``surface area'' to which abductive reasoning applies is comparatively small and carefully delineated.
Once the reader is confident that the definitions align with their internal mental model, which can be done with abductive reasoning as discussed above, they may choose how far to dig into the associated proofs while remaining confident that the ITP ensures correctness.
Overall, this has the effect of \emph{separating} the exposition of the material from the issue of trusting its correctness.

\section{Conclusion}

In \S\ref{section:sources-of-complexity} we discussed two kinds of sources of complexity in mathematics: \emph{accidental} and \emph{inherent}.
We also saw that ITPs can manage and remove accidental complexity, essentially \emph{by default}.
On the other hand, in \S\ref{section:spec-driven-dev} we saw that spec driven development and abstraction boundaries, in conjunction with an ITP, can also help in taming \emph{inherent} complexity.
As discussed in that section, a fundamental feature of spec driven development is the concept of \emph{mathematical debt} which can be effectively tracked using an ITP.
While this debt has to be repaid at some point, it is fundamentally \emph{non-blocking} for the project as a whole.

The targets arising in spec driven development, along with other mathematical objects in question, are frequently interconnected in highly nontrivial ways.
With the use of an ITP, contributors only need to keep track of the objects which are \emph{locally} relevant, thereby reducing the cognitive load significantly.
At the same time, these targets can be recursively broken into smaller components until they can be handled directly.
This is particularly useful in collaborations as this allows for the work to be easily distributed among contributors with potentially varying areas of expertise.

Another key property of spec driven development is that reliance on abstraction boundaries becomes effectively \emph{required} with this approach.
As a byproduct, refactors generally become much easier.
As the ITP ensures consistency within the project at all times, it also keeps track of the changes that must be made as a consequence of such refactors.

We feel that there is a lot of potential to develop tools to help facilitate spec driven development in pure mathematics, in conjunction with an ITP.
Above we discussed how an ITP ``keeps track'' of the complicated dependency graph of objects (including targets) within a mathematical project.
Right now, this is essentially done in the process of \emph{type-checking}.
However, a tool that does this in a more \emph{explicit} manner, by identifying precisely the various objects associated with a given piece of mathematics within an ITP, would be even more beneficial as it would tell contributors \emph{precisely} what they need to focus on when working on individual targets.
Patrick Massot's blueprint software~\cite{blueprint-software}, which was used in LTE,%
\footnote{\url{https://leanprover-community.github.io/liquid/}}
can be seen an initial approximation to this idea.

Another tool that would be useful in projects such as LTE is a system that would keep track of the ``weights'' of targets.
Essentially, we envision that each \sorry~in a given project would be tagged with a ``weight'' that is meant to act as a quantitative approximation for its complexity.
Progress in the project as a whole can then be tracked by computing the total remaining weight.
As the history of projects such as LTE is stored using a version control system (e.g.~\texttt{git}), the evolution of the weights of targets (along with other metadata) could be mined and studied, and potentially used in machine learning applications.

In this article, we have painted a very positive picture about the current use of ITPs in pure mathematics.
The formal verification of significant results, such as Hales' proof of the Kepler conjecture~\cite{zbMATH06750549} and the liquid tensor experiment, have clearly demonstrated that ITPs are valuable tools for the verification of hard proofs.
To summarize, we discussed how the use of abstraction boundaries and spec driven development within the rigid framework of an ITP helps
\begin{itemize}
  \item facilitate collaboration amongst mathematicians with different expertise;
  \item manage accidental complexity and thus reduce cognitive load;
  \item provide a guarantee of formal correctness of proofs, thereby reducing the amount of material that needs to be trusted.
\end{itemize}
We envision that these methods can also provide similar benefits in original research.

However, we acknowledge that there are still a number of high costs that come with the use of ITPs and related tools that pose a barrier to wider adoption within mathematics.
For example, they currently have a steep learning curve, and the difference between an informal piece of mathematical exposition and its formalization is still too large.
Nevertheless, ITPs and their formalized mathematics libraries have seen tremendous improvements in the last few years, while the successes in the field have helped accelerate progress.
We are also optimistic that advances in metaprogramming, artificial intelligence, and type theory itself, will help resolve some of the remaining issues in the near future.

\bibliographystyle{amsplain}
\bibliography{refs}
\end{document}